\documentclass{nm}
\renewcommand\thisnumber{x}
\renewcommand\thisyear {201x}

\renewcommand\thisvolume{xx}
\usepackage{charter}
\usepackage[charter]{mathdesign}
\usepackage{amsmath}
\usepackage{amsfonts}
\usepackage{graphicx}%
\usepackage{latexsym ,rawfonts}
 \usepackage{multirow}
\usepackage{tikz}
\usetikzlibrary{shapes}
\usepackage{color}

\usepackage{algorithmic}
\usepackage{algorithm}
\usepackage{subfigure}
\usepackage{float}
\usepackage{geometry}                
\usepackage{graphicx}
\usepackage{algorithmic}
\usepackage{algorithm}
\usepackage{enumerate}
\newcommand{\wm}{\omega_{max}}
\newcommand{\be}{\begin{equation}}
\newcommand{\ee}{\end{equation}}
\newcommand{\thf}{\theta_1^h}
\newcommand{\ths}{\theta_2^h}
\newcommand{\tthf}{\theta_1^H}
\newcommand{\tths}{\theta_2^H}
\newcommand{\bw}{{\bf \omega}}
\newcommand{\bbw}{|{\bf \omega}|}
\newcommand{\wo}{\omega_1}
\newcommand{\wt}{\omega_2}

\newcommand{\tahh}{\tilde{L}^H}
\newcommand{\tmhh}{\tilde{M}^H}
\newcommand{\lh}{L^H}
\usepackage{graphicx}
\usepackage{subfig}

\begin{document}
\markboth{ I. Livshits}{Multigrid for Helmholtz}
\title{Shifted Laplacian based multigrid preconditioners  for solving  indefinite Helmholtz equations}
\author[Ira Livshits]{Ira Livshits\affil{1}}
\address{\affilnum{1}\ Department of Mathematical Sciences, Ball State University, Muncie IN, 47306, USA.
       }

\emails{
 {\tt ilivshits@bsu.edu} (lra Livshits)
 }

\begin{abstract}
Shifted Laplacian multigrid preconditioner \cite{vuik2} has become a tool du jour for solving highly indefinite Helmholtz equations. The   idea is to  add 
a complex damping  to the original  Helmholtz operator and  then apply a multigrid processing to the  resulting operator using it   to precondition  Krylov  methods, usually  Bi-CGSTAB.  Not only such preconditioning  accelerates  Krylov iterations, but  it does so  more efficiently  than the  multigrid applied to original Helmholtz equations. In this paper, we  compare  properties of   the Helmholtz operator with  and without the shift and propose  a new combination of the two.  Also applied  here  is a relaxation of normal equations that replaces diverging  linear schemes on some intermediate scales. 
Finally,  an  acceleration  by the ray correction \cite{BL97} is considered. 
\end{abstract}




\maketitle
\section{Introduction}
Considered here is a two-dimensional Helmholtz equation
\be
Lu = \Delta u(x)  + k^2(x) u(x) = f(x), \quad x\in \Omega \subset \mathbf{R}^2,
\label{eq:diff}
\ee
accompanied by the first-order Sommerfeld boundary conditions
\be
\frac{\partial{u(x)}}{\partial{n}}-iku(x) = 0, \quad x \in \partial \Omega,
\ee
where $n$ is an outward normal.
 Discretized on 
 a sufficiently fine scale $h$,   $ k h \le 2\pi/10$,  using  standard discretization methods,   (\ref{eq:diff}) yields a system of linear equations
\be
L^h  u^h = f^h,
\label{eq:disc}
\ee 
where $L^h \in \mathbf{C}^{N\times N}$ is a sparse matrix, where   $N$  is  typically very  large.

Different methodologies  applied to  (\ref{eq:disc}) range from direct, e.g. \cite{rokhlin, eng1} to iterative ones, including multigrid.  The latter often offers a  high approximation accuracy at optimal computational costs.  Multigrid approaches  for  (\ref{eq:diff})  notably include \cite{ Elman, vuik2, haber11, barry,vanek2} among others. The most practical  multigrid  method to date  is the Shifted Laplacian approach e.g.,  \cite{vuik1, vuik2}. It employs  a discretization  of a shifted differential operator $M = L + i k^2\beta$,
\be
M^h = L^h + ik^2\beta,
\label{eq:mh}
\ee
as a preconditioner to $L^h$,  with $i = \sqrt{-1}$ and typical $\beta = 0.5$ as assumed  throughout the paper.  The complex damping 
 helps with some of the challenges presented by the Helmholtz operator, it is easy to implement, and, most importantly,  $M^h$ based multigrid preconditioner  
 significantly accelerates  Krylov  iterations. 
Another obvious idea,  justly  overlooked  due its poor performance, is  applying   multigrid directly to (\ref{eq:disc}).  In this paper the two approaches, based on  the Helmholtz  and the Shifted Laplacian operators, are compared, and  a hybrid method  is proposed.   Also briefly discussed is  the ray correction  \cite{BL97}.

 The remainder of the paper is organized as follows.  Operator (\ref{eq:diff})  and  error components, whose treatment is  essential to effectively  solving it, are discussed in Section \ref{sec:errors}. The Helmholtz (HLM) and the Shifted  Laplacian (SL) approaches  are  compared from two perspectives: how accurately  $L^H$ and $M^H$, $H = 2h, 4h, \dots$ approximate the finest grid operator $L^h$,  Section \ref{sec:symbols},  and how well Gauss-Seidel relaxation, applied to $L^H$ and $M^H$, $H = h, 2h, 4h, \dots$  converges for different types of error components,  Section \ref{sec:smoothing}.  An optimal strategy which involves combining the two methods is  suggested in Section \ref{sec:optimal};  numerical  experiments are presented  in   Section \ref{sec:numerics}, and the  concluding remarks are given   in Section \ref{sec:conclusions}.

\section{Error components and the Helmholtz operator}
\label{sec:errors}

Any efficient  multigrid algorithm works in the following way:  each   coarse grid  operator $A^H$, $H=2h, 4h, \dots$     approximates  the finest grid operator   $A^h$ for all components unreduced by processing on finer grids;   error  $e^H$ with large relative residual \be \| A^H e^H\| \gg  \|e^H\|\ee is  practically annihilated  by a few relaxation sweeps applied to 
\be A^H e^H = r^H, \ee
where $r^H$ is the coarse grid residual, an average of the residual computed on the finer scale, $H/2$. The remaining error, with small relative residual,  is   accurately approximated on the next  coarser scale, $2H$, and so forth.  This means in particular that  error components with the smallest relative residuals,  i.e.,  the {\sl near-kernel error components} of $A^h$,
\be A^h e^h \approx 0,
\label{eq:nkc}
\ee
have to be approximated on all, including the coarsest, scales, which  works naturally when they  
are smooth. This is not the case Helmholtz operators with large wave numbers.  There components (\ref{eq:nkc}) are of the form (at the interior) 
\be e(x,y) = e^{i(\wo x + \wt y)}, \quad   
\label{nkc}
\ee
with $|\bw |  = \sqrt{\wo^2 + \wt^2} \approx k$.  (In further discussion,  instead of a  general  $\bbw \approx k$, a more specific   $(1-\alpha_0) k  \le |\bw| \le (1+\alpha_1) k$,  for some  $0 < \alpha_0, \alpha_1 < 1$, is used.)  Starting with some scale $H$,  these components become oscillatory; for larger $k$ it happens  on finer  $H$. 

Next  properties  of $L^H$  and $M^H$ when applied to different  error components are analyzed  and compared.

\section{Approximation by Helmholtz and  Shifted Laplacian operators}
\label{sec:symbols}
An  approximation accuracy of a fine-grid operator by a coarse-grid operator  is often measured by comparing {\sl symbols} of the two for Fourier components visible on the coarser scale.  Generally, a   symbol of an  operator $A$ applied to  $e^{i(\wo x + \wt y)}$  is defined as a complex coefficient  $\tilde{A}(\wo, \wt)$:
\be A e^{i(\wo x + \wt  y)}  \approx \tilde{A}(\wo, \wt)  e^{i(\wo x + \wt  y)}. \ee

\begin{figure}
 \includegraphics[width = 5in]{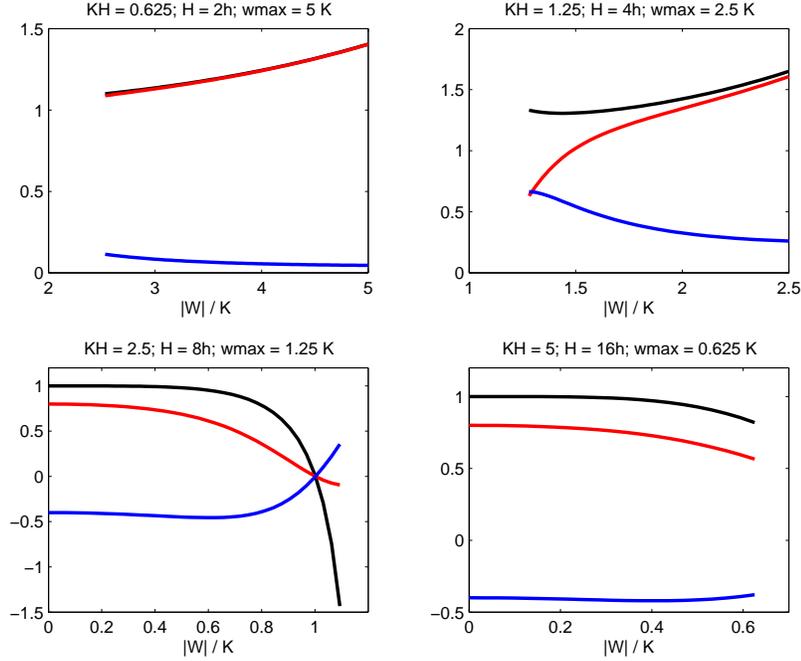}
 \caption{Symbol ratios  for Fourier components visible on scale  $H$;  the $x$-axis  variable is   $\bbw/k$.  The black (top) line  shows the HML  ratios (\ref{symb:helmholtz}); the red line (middle)  and the blue line (bottom) are  the real and the imaginary parts of the SL symbol ratios (\ref{symb:shifted}). The ratios shown are for components with $\wo = \wt$.}
 \label{fig:symbols}
\end{figure}

\noindent
For a coarse-grid  correction either by $L^H$ or by $M^H$   to provide an adequate approximation to solution of (\ref{eq:disc}) the symbol ratios, defined with  $\thf = \omega_1 h, \ths = \omega_2 h, \tthf = \omega_1 H $ and $\tths = \omega_2 H$, 
 \be
\tau^{HLM}_{H}(\wo, \wt) = \frac{\tilde{L}^h(\wo, \wt)}{\tahh(\wo, \wt)} =\frac{(2\cos \, \thf+2\cos \, \ths - 4 +k^2h^2)\, H^2}{(2\cos \,\tthf+2\cos \,\tths - 4 +k^2H^2) \, h^2}
 \label{symb:helmholtz}
 \ee
 and 
\be
 \tau^{SL}_{H}(\wo, \wt) = \frac{\tilde{L}^h(\wo, \wt)}{\tmhh(\wo, \wt)}  \frac{(2\cos \, \thf+2\cos \, \ths - 4 +k^2h^2)\, H^2 }{(2\cos\, \tthf+2\cos \,\tths - 4 +k^2H^2(1+i\beta))\,h^2}
  \label{symb:shifted}
 \ee
should be close  to one. 
\noindent
To illustrate how values of (\ref{symb:helmholtz}) and (\ref{symb:shifted}) change when considered on  increasingly coarser  scales,   Figure \ref{fig:symbols}  
shows  results  for   error components,  that are oscillatory on each  scale $H$, 
  $ \pi/2\le \wm H \le \pi$,  where $\wm = \max\{\wo, \wt\}$. The exception is the last  subfigure  which    depicts the entire range  visible on $H = 16h$,   $ 0 \le \wm H \le \pi$.

As Figure \ref{fig:symbols}  suggests, operators $L^H$ and $M^H$  exhibit    similar accuracy for  high-frequency components but  differ for the near-kernel   (\ref{nkc}) and for  lower frequencies.  More precisely, for 
\begin{itemize}
\item $kH \le  0.625$: all components with $\pi/2 <  \wm H \le \pi $  are well approximated by  $L^H$ and $M^H$:   $\mathop{\rm Re}(\tau_H^{SL}) \approx \tau_H^{HLM} $ are    close to one,  $\mathop{\rm Im} (\tau_H^{SL}) \le .1$; 
\item $kH = 1.25$: all components with $\pi/2 <  \wm H \le \pi $ 
  are accurately   approximated by $\lh$, they satisfy  $(1+\alpha_1) k \le |\bw| \le 2(1+\alpha_1) k$, $\alpha_1 \approx .4$. The accuracy deteriorates for smoother components, in particular   as $|w|$ approaches  $k$.  $M^H$ provides an accurate approximation for a smaller range of components,  the ones with $(1+\beta_1)\le |\bw|  \le 2(1+\alpha_1) k$, $\beta_1 \approx .8$.  The growing imaginary part of $\tau_H^{SL}$ for smaller $|\bw|$ affects  the  approximation quality.  
 Both $L^H$ and $M^H$  fail to approximate  components (\ref{nkc}) though in a different way.\footnote{The wrong approximation and relaxation  of these components by the SL operators is  an asset when the SL approach is  used as a preconditioner, as it regroups  the eigenvalues corresponding to such components in a way that makes them  more treatable by Krylov  methods \cite{vuik2}. }
\item $kH = 2.5$: $L^H$ provides an accurate approximation for $|\bw|  \le (1-\alpha_ 0) k$, with $\alpha_0 \approx .1 $ and does not  approximate components with $|\bw| \ge k$;  $M^H$ gives a rise to a wrong approximation for  all components in question, though manages to do so in the right way (see the footnote); 
\item $kH = 5$:  all components visible on scale $H$, $0 \le |\bw| \le 0.65  k$,   have an accurate  approximation by $L^H$,  but not by $M^H$  (due to a large negative imaginary part).
\end{itemize}

To summarize,    a sequence of  coarse-grid Helmholtz operators $\{L^H\}_{H>h}$ accurately approximates the finest grid Helmholtz operator $L^h$ for  all Fourier components except (\ref{nkc}), more precisely with $(1-\alpha_0)  k \le |\bw| \le (1+\alpha_1) k$, with $\alpha_0 \approx .1$ and $\alpha_1 \approx .4$.
 Coarse-grid Shifted Laplacian  operators $\{M^H\}_{H>h}$ approximate the finest-grid Helmholtz operator  for all oscillatory components, failing to approximate  both  (\ref{nkc})  and (unlike $L^H$)  smooth error,  more precisely, components  with $ 0 \le |\bw| \le (1+\beta_1) k$, with  $\beta_1 \approx .8$.

\section{Gauss Seidel relaxation for $L^H$ and $M^H$}
\label{sec:smoothing}
\begin{figure}
  \includegraphics[width = 5in]{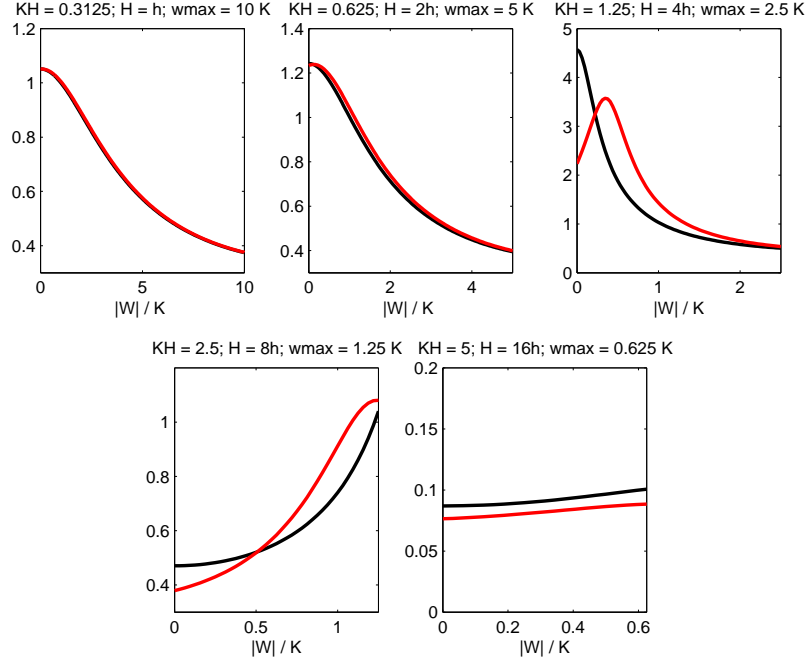}
  \caption{The black line shows the rates for Helmholtz operators; the red line - for  Shifted  Laplacian ones;
          convergence rates of  Fourier components visible on scale  $H$; the finest scale satisfies $kh = 0.3125$; 
      the $x$-axis variable is $|\bw|/k$; results are shown for Fourier components with $\wo = \wt$; 
       }
       \label{fig:smoothing}
 \end{figure}

Application of one iteration of the lexicographic Gauss-Seidel relaxation  to $\lh$ and  $M^H$ yields the following amplitude change of  an erroneous  Fourier component $e^{i(\wo x + \wt y)}$ 
\be \mu_H^{HLM}(\tthf, \tths) = \biggl| \frac{\exp(-i \tthf) + \exp(-i \tths)}{\exp(i \tthf)+\exp(i\tths)-4+k^2 H^2}\biggr |
\label{conv:helmholtz}
\ee
and  
\be \mu_H^{SL}(\tthf, \tths) = \biggl| \frac{\exp(-i \tthf) + \exp(-i \tths)}{\exp(i \tthf)+\exp(i\tths)-4+k^2(1+i \beta) H^2}\biggr |
\label{conv:shifted}
\ee
for Helmholtz and  Shifted Laplacian operators, respectively.
Typically, in predicting  a convergence rate  of a multigrid solver,  the {\sl smoothing}  properties of the relaxation is the main parameter. It is measured by a {\sl smoothing} factor: 
\be 
\tilde{\mu}_H = \max_{\pi/2 \le \max\{|\tthf|, |\tths|\} \le \pi}  \mu_H(\tthf, \tths).
\ee
\noindent
For  $\lh$ and $M^H$ there is an additional phenomenon --  divergence of  smooth error component.  To monitor that,    an overall  convergence  rate is also considered:
\be\hat{\mu}_H = \max_{0 \le \max\{|\tthf|, |\tths|\} \le \pi}  \mu_H(\tthf, \tths);   \ee
  $\hat{\mu}_H > 1$ means divergence. Figure \ref{fig:smoothing} shows $\mu_H(\wo, \wt)$ for $L^H$ and $M^H$ on increasingly coarse scales starting with the finest,  $kh = 0.3125$.  It suggests that Gauss-Seidel relaxation   performs similarly when applied to $L^H$ and $M^H$. In particular for 
\begin{itemize}
\item  $kH \le 0.3125$: 
$\tilde{\mu}_{H}^{HLM} \approx \tilde{\mu}_{H}^{SL}   \approx  0.5$ and 
$ \hat{\mu}^{HLM}_H \approx \hat{\mu}_H^{SL} \le  1.05$;  
\item  $kH =  0.625$: 
 $\tilde{\mu}_{H}^{HLM} \approx \tilde{\mu}_{H}^{SL}   \approx  0.7$ and  
$ \hat{\mu}^{HLM}_H \approx \hat{\mu}_H^{SL} \le  1.1$;  
\item   $kH = 1.25$:  divergence of  smooth error components becomes prohibitively large,  with  $\hat{\mu}^{HLM}_H \approx 4.5$  and $\hat{\mu}^{SL}_H \approx 3.5$;  no error reduction for $|\bw| \approx k$.  However, error components 
with $|\bw | \ge 1.3 k$ for $L^H$ and with $|\bw| \ge 1.8k$ for $M^H$ are reduced by at least  the factor of  $0.7$; 
\item    $kH = 2.5$:  $\tilde{\mu}^{HLM} \approx \tilde{\mu}^{SL} \approx 1$ -- no  convergence  for  (\ref{nkc}); for smooth  components, $\mu^{HLM}_H(\wo,\wt) \le 0.7$ for    ($|\bw | \le .9 k$) and  $\mu^{SL}_H(\wo,\wt) \le 0.7$  for ($|\bw | \le .8 k$).
\item   $kH = 5$: $\hat{\mu}^{HLM}_H \approx .1 $ and  $\hat{\mu}^{SL}_H \approx  .085 $ making  a  few relaxation sweeps  an equivalent to a direct solver; no coarser grids are needed. \\
\end{itemize}

\noindent
Overall,  Gauss-Seidel relaxation for both approaches performs well on scales  with $kH \lessapprox 0.625$ and $kH \gtrapprox 2.5$. 
It fails to reduce near-kernel components (\ref{eq:nkc})  on any grid  and  diverges  smooth error components  when $kH \approx 1.25$. To avoid or diminish the latter effect,  Gauss-Seidel is applied to the normal operator 
$(L^H)^T L^H$  or $(M^H)^TM^H)$ instead of original $L^H$  or $M^H$, where $^T$ here means transposed, complex conjugate.  This is done in the spirit of  Kaczmarz iterations \cite{kaczmarz}  known to be slow but convergent.  The number of relaxation sweeps  on  this  scale is  higher than on others.
 
 \begin{remark}
 The actual constants in the discussion above as well as in Section \ref{sec:symbols} are partial for the chosen parameters; they aim at giving  a qualitative 
 understanding of the processes described. While the study is conducted for Gauss-Seidel iterations, similar conclusions, with slightly different constants, can be made for other linear iterative schemes such as Jacobi or SOR.\\
 \end{remark}

\section{Optimal algorithm}
\label{sec:optimal}
A multigrid  V-cycle is applied to (\ref{eq:diff})  in  three variants. It employs: 
\begin{itemize}
\item Operators $\lh$ and/or  $M^H$,  second-order FD discretizations of $L$ and $M$  with   five-point stencils; 
\item  bilinear interpolation;
\item full weighting;
\item Gauss-Seidel iterations:
\begin{itemize}
\item  one pre- and post-smoothing steps  on all scales except  $kH \approx 1.25$, applied either to $L^H$ or to $M^H$; 
\item  four   pre-  and  post-smoothing  steps on scale  $kH \approx 1.25$, applied either to $(L^H)^T L^H$ or to $(M^H)^T M^H$. 
\end{itemize}
\end{itemize}
On each scale a coarse-grid operator is  used in two capacities: 
\begin{enumerate}[(A)]
\item for relaxation; 
\item for  computing  coarse-grid residuals.
\end{enumerate}
\noindent
Three variants are considered: 
\begin{itemize}
\item  HLM-V employs  $L^H$ both for (A) and for (B); 
\item  SL-V employs  $M^H$   both for (A) and for (B)
\item  HYB-V always employs $L^H$ for (B). $L^H$ is also used  for (A) on all grids  except $0.625 \lessapprox   kH \lessapprox 1.25 $ where  it is replaced by $M^H$. 
\end{itemize}
The motivation for the hybrid method comes from observations reported in Sections \ref{sec:symbols} and \ref{sec:smoothing}  concerning   performance of  SL and  HLM operators on intermediate and coarse scales. (On finer grids both act very similarly, and either one can be used.) The strength of the Shifted Laplacian approach, studied in detail in \cite{vuik1,vuik2}, is the transformation (not reduction) of the near-kernel error components,  that mostly occurs  on intermediate scales. This is the reason  for employing  $M^H$ in relaxation there. 

On coarse grids, however,  Helmholtz operators $L^H$, give  a rise to an accurate approximation of smooth components, and, together with a fast convergence  by  Gauss-Seidel there, allow for an  efficient coarse-grid  correction.  Therefore, $L^H$ is used in relaxation on the coarsest scale(s).

\section{Numerical Experiments and Computational Costs}
\label{sec:numerics}
The V-cycle based  variants, along with   the original Shifted Laplacian  (OSL) \cite{vuik2} multigrid preconditioner, are compared, and   their   computational costs are discussed. Bi-CGSTAB serves as an outer iteration.
 Also briefly introduced is the idea of the ray correction \cite{BL97}, and numerical results for   HLM,  SL and  HYB, enriched by it, are presented.  
\subsection{Numerical results}
First, the algorithms are tested for (\ref{eq:diff}) with a constant $k$, considered on $\Omega = [0,1]^2$, and the  results are presented in  Table \ref{tab1}.  Initial approximations $x^0$ are zero in all experiments;  iterations are performed until the initial residual $\|r^0\| = \| f\|$ is reduced by a factor of $10^{7}$.
In Tables \ref{tab1}-\ref{tab2},   right-hand-sides are  homogeneous except at the center of $\Omega$, where  $f(.5,.5)=1$;
\begin{table}[!h]
\begin{center}
\begin{tabular}{l |c c c c c}
$k$ & 40 & 50 & 80 & 100 & 150 \\ \hline
$h$ & $1/64$ & $1/80$ & $1/128$ & $1/160$ & $ 1/240$ \\ \hline
OSL  & 26 & 31 & 44 & 52 & 73 \\ \hline
SL-V & 19 & 24 & 27.5 & 31  & 38 \\ \hline
HYB-V & 16 &  20.5 & 23 &  26.5 &  31.5 \\ \hline
\end{tabular}
\end{center}
\caption{The number of  Bi-CGSTAB iterations for different preconditioners and values of  constant wave numbers; in all experiments $kh = 0.625$. 
}
\label{tab1}
\end{table}

The results show that both the  SL-V and HYB-V  preconditioners  are more  efficient   than OSL,  and  the hybrid approach outperforms the Shifted Laplacian.

In Table \ref{tab2},  performance of  SL-V and  HYB-V methods is tested for  the same model problem when  considered on  increasingly finer $h$; both  show  an  improved  convergence while computing   increasingly accurate solutions. \\

\begin{table}[!htb]
\begin{center}
\begin{tabular}{l | c c c r}
$h$ & 1/64 & 1/128 & 1/256 & 1/512 \\ \hline
SL-V&   19 & 18  & 17.5    &  16 \\ \hline
HYB-V  &  16  &  15.5  &  15  &  14 \\ \hline
\end{tabular}
\end{center}
\caption{The number of Bi-CGSTAB iterations;  $k=40$, $kh$  ranges from $0.625$ to $0.078125$.}
\label{tab2}
\end{table}

\begin{figure}[!htb]
 \centering
 \includegraphics[scale=0.6]{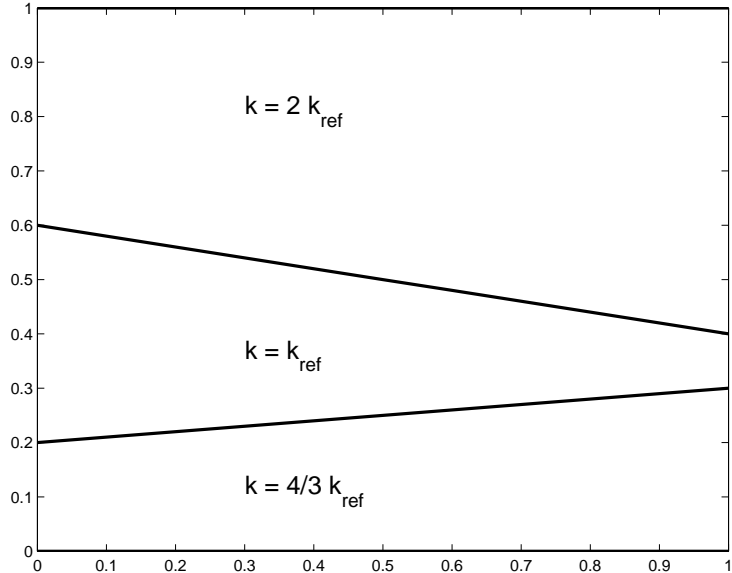}
 \caption{Wave number distribution in the wedge problem The point source is located at the middle of the upper boundary: $f(0.5,1.) = 1$ and zero elsewhere. The choice of  $k(x)$ is in spirit of the wedge example in \cite{vuik2}. The domain remains $[0,1]^2$.
 } \label{fig5}
 \end{figure}

\noindent
Next  considered is (\ref{eq:diff})  with a  heterogeneous medium - a   wedge problem shown  in Figure \ref{fig5}, with  numerical experiments  presented in Table \ref{tab6}. Again,  the hybrid   preconditioner  performs better  than the  Shifted Laplacian  does.

 \begin{table}[!htb]
 \begin{center}
 \begin{tabular}{l | c c c c r}
 $k_{ref}$ & 15  & 30 & 60  & 120  & 240 \\ \hline
 SL-V &  13 & 18.5 &  33  &  49.5   &   61 \\ \hline
 HYB-V  & 9.5 &  14  &   23  &  36.5  &  41  \\ \hline
 \end{tabular}
 \end{center}
 \caption{The number of  Bi-CGSTAB iterations for  SL-V and  HYB-V preconditioners;  in all experiments $ k_{ref} h \approx 0.2344$; the value of $k$ in the Table varies from $15$  to  $480$. }
 \label{tab6}
 \end{table}

Our experiments are performed for  a slightly different problem that the ones reported  in \cite{vuik2}.  We mention, however, that in \cite{vuik2} the experiments were performed for $k \le k_{max} =  240$, which arises for $k_{ref} = 120$,  and  it required $66$  Bi-CGSTAB iterations. Our experiments with the same  $k_{ref} = 120$ (and  $k_{max} = 240$)  require  only  $36.5$ Bi-CGSTAB iterations with HYB preconditioner.  
 
Noticeably missing from  action so far is  HLM-V approach,  and this is because its acceleration of  Bi-CGSTAB or other Krylov methods, is inferior to the SL-based algorithms. This changes, however, when the ray multigrid approach \cite{BL97} is used as an additional  coarse-grid correction,  \cite{Thesis, BL97, BL06}.  It is  based on the assumption that the near-kernel error components  (\ref{nkc}) can be represented as 
\be e = \sum_{j=1}^8 \hat{a}_j e^{i(k_1^j x + k_2^j y)}, \quad  (k_j^1, k_j^2) = k\biggl(\cos \,  \frac{j\pi}{4}, \sin \, \frac{j\pi}{4}\biggr),
\ee 
with smooth  {\sl ray}  functions $\hat{a}_j$. The idea is than to reduce the task of  computing  $e$ to a much easier task of approximating  each  $\hat{a}_j$ individually on some coarse scale.  This process itself   reduces a range of the near-kernel Fourier error components with $(1-\gamma_0) k \le |\bw| \le (1+\gamma_1) k$. Constants $\gamma_0, \gamma_1$  depend on  relaxation strategy and problem parameters:  typical values  are  $\gamma_0 \approx \gamma_1 \approx 0.3$. This means that all error components not well approximated/well reduced by  HLM-V are in this range, and   they are all treated  by the ray correction.
Results for HLM-V,  HYB-V  and  SL-V cycles,  accelerated  by the ray correction, are presented in Table \ref{tab7}.  No Krylov outer iterations  are employed: each method   serves as a solver  rather than  a preconditioner;  HLM-V cycle  with the ray correction is the original wave-ray algorithm. The cost of each iteration in this Table is about twice lower than  iteration costs in other Tables, where one Bi-CGSTAB employs a multigrid preconditioner {\sl twice}.
\begin{table}[!htb]
 \begin{center}
 \begin{tabular}{l | c c c  r}
 $k$ & 20 & 40  & 80 & 160 \\ \hline
HLM-V  & 16 & 16  & 17 & 18 \\ \hline
 SL-V &   23   & 34& 41  & 48 \\ \hline
 HYB-V  & 26 &  31 & 37 & 39 \\ \hline
 \end{tabular}
 \end{center}
 \caption{The number of V-cycles using  HLM-V,  SL-V  and  HYB-V approaches enhanced by the ray correction;   in all experiments $ kh = 0.3125$; $f(.5,.5) =0$, and it is zero elsewhere. }
 \label{tab7}
 \end{table}

\subsection{Computational Costs}
\label{sec:costs}
Costs of  SL-V,  HLM-V, and HYB-V preconditioners are close to costs of a standard multigrid  $V(1,1)$ cycle applied to  a Laplace operator; the main  difference is the cost of  the extra six relaxation sweeps applied to the normal equation on scale with $kH \approx 1.25$ (Six is eight per level minus standard two per grid in $V(1,1))$. While the absolute cost of these iterations remains the same for a given (\ref{eq:diff}), its relative fraction in the overall costs  becomes smaller when  (\ref{eq:disc}) is discretized on finer scale $h$. The OSL preconditioner is implemented differently from the algorithms  discussed here: it employs a $F(1,1)$ cycle in the algebraic multigrid framework using the operator dependent-interpolation based on de Zeewv's transfer operators \cite{zee}.
The $F(1,1)$  cycle becomes more expensive (in computational costs) than our  almost $V(1,1)$ cycle starting with $kh = 0.3125$ and finer.

\section{Conclusions} 
\label{sec:conclusions}
Standard multigrid V-cycle is applied to the Helmholtz and the Shifted Laplacian operators, and the resulting algorithms are employed  as  preconditioners for   Bi-CGSTAB, used to solve the indefinite Helmholtz equations. The Shifted Laplacian approach shows a superior performance. However, after  analyzing   approximation and  relaxation properties of both operators, a hybrid method, a combination of the two, is proposed, yielding an improved convergence.  With the ray correction, the HLM approach works significantly better - resulting in a well  scalable algorithm with convergence nearly independent on  wave numbers. 
\bibliographystyle{plain}
\bibliography{helm2.bib}	
\end{document}